\documentclass[12pt]{article}
\title{Harmonic analysis of finite lamplighter random walks.}

\author{Fabio Scarabotti, Filippo Tolli}

\usepackage{latexsym,amsbsy,amsmath,amscd,amssymb}

\usepackage{amsthm}
\usepackage{amssymb,amsmath,latexsym}
 \newtheorem{definition}{Definition} [section]
 \newtheorem{remark}[definition]{Remark}

 \newtheorem{proposition}[definition]{Proposition}
 \newtheorem{theorem}[definition]{Theorem}
 \newtheorem{corollary}[definition]{Corollary}
  \newtheorem{lemma}[definition]{Lemma}

\begin{document}

\maketitle

\begin{abstract}
Recently, several papers have been devoted to the analysis of
lamplighter random walks, in particular when the underlying graph is
the infinite path $\mathbb{Z}$. In the present paper, we develop a
spectral analysis for lamplighter random walks on finite graphs. In
the general case, we use the $C_2$-symmetry to reduce the spectral
computations to a series of eigenvalue problems on the underlying
graph. In the case the graph has a transitive isometry group $G$, we
also describe the spectral analysis in terms of the representation
theory of the wreath product $C_2\wr G$. We apply our theory to the
lamplighter random walks on the complete graph and on the discrete
circle. These examples were already studied by Haggstrom and
Jonasson by probabilistic methods.
\footnote{{\it AMS 2002 Math. Subj. Class.}: Primary: 43A85; secondary: 05C05, 20C15, 20E22, 60G50.}
\end{abstract}

\section{Introduction.}

Let $X$ be a simple, locally finite, connected graph. Put in each
vertex a lamp, which may be on or off. A lamplighter performs the
simple random walk on $X$ and when he moves from a vertex $x$ to a
vertex $y$ he changes randomly the state of the lamps in $x$ and $y$,
that is both the lamps may be turned on or off with equal
probability. In other  words, we may construct a new graph whose
vertex set is $\mathcal{L}(X)=\{(\theta,x)|\theta:X\rightarrow
\{0,1\},x\in X\}\equiv C_2^X\times X$ and two vertices $(\theta,x)$
and $(\sigma,y)$ are connected if $x$ and $y$ are connected in $X$
and $\theta\equiv\sigma$ in $X\setminus\{x,y\}$.
Then $\theta(x)$ is the state of the lamp in $x$, $\theta(x)=0$ if it is off, $\theta(x)=1$ if it is on, and the lamplighter random walk on $X$ is just the simple random walk on $\mathcal{L}(X)$.\\

These kinds of processes (that have many variants) have been studied by many authors; in particular, we mention \cite{Ba-Woe,Di-Sc,Gri-Zu} devoted to the spectral analysis of the lamplighter random walk on the infinite path; actually, the paper of L. Bartholdi and W. Woess treats the more general case of the Distel-Leader product of two infinite homogeneous trees (see also \cite{Woess2}).
We also refer to \cite{Woess}, that contains a more general construction where the lamps are replaced by the vertices of another graph. On the other hands, the finite case has been treated by O. Haggstrom and J. Jonasson \cite{Ha-Jo}, who analyzed by probabilistic techniques the lamplighter processes on the complete graph and on the discrete circle,  and by Y. Peres and D. Revelle \cite{Pe-Re}, who used analytic techniques for the lamplighter process on finite tori.\\

In the present paper, we develop a suitable spectral analysis for lamplighter random walks on finite graphs. We start from the following simple observation: the lamplighter random walk is $C_2^X$-invariant. This group acts on the lamps coordinatewise: if $\theta\in C_2^X$ and $(\omega,x)\in\mathcal{L}(X)$ then $(\theta + \omega)(x)=\theta(x)+\omega(x)$ $\text{mod}\;2$, and the action on $\mathcal{L}(X)$ is simply

\begin{equation}\label{gensym}
\theta\cdot(\omega,x)=(\theta+\omega,x).
\end{equation}

Clearly, this is not a transitive action. In \cite{Sc-To2}, starting from the results in \cite{BDPX}, we developed a suitable harmonic analysis for a finite Markov chain with a nontransitive group of symmetries. In the present setting, our methods simplify noticeably: when we restrict the Markov operator to the isotypic components of the permutation representation of $C_2^X$ on the lamplighter graph $\mathcal{L}(X)$, we get a series of eigenvalue problems on the graph $X$, that lead to a complete spectral analysis of the lamplighter chain. \\

The plan of the paper is the following. In Section 2 we establish a series of notation used in the paper. In Section 3 we analyze the lamplighter process described above. In Section 4 we analyze the lamplighter process on the complete graph on $n$ vertices. In particular, using the standard techniques developed by P. Diaconis \cite{Diaconis}, we show that the chain has a cut-off after $k=\frac{1}{2}n\log n$ steps. This result has already been obtained in \cite{Ha-Jo} by means of purely probabilistic techniques. In Section 5 we analyze one of the possible variations of the lamplighter construction:  we put the lamps on the edges of the graph. When we move form $x$ to $y$, the lamp in the edge $\{x,y\}$ is randomized. In Section 6 we compute the spectrum of the lamplighter random walk on the discrete circle, with the lamps on the edges. The result in this section may be considered as a finite analogous of the computations on the infinite path; moreover, both the finite and infinite cases are random walks on groups, namely the wreath products $C_2\wr C_n$ and $C_2\wr \mathbb{Z}$. The spectral computations in this section have a clear connection with the eigenvalue problems on finite trees treated in \cite{Strang2,Sca,Sc-To}.  We have written Sections 3-6 with a minimum of group formalism, in order to make this part of the paper accessible with only a discrete/probabilistic background. In the remaining part of the paper, we make a systematic use of group representation theory. In Section 7, we prove a general decomposition theorem for the permutation representation of a group $G$ on a space of the form $C_2^Z\times X$, where both $X$ and $Z$ are $G$-homogeneous spaces. This is more than is needed for the lamplighter random walks; in fact, in Section 9 we show that a decomposition derived by Schoolfield (for the Bernoulli-Laplace diffusion model with sign) may be easily deduced from our general result. In Section 8, we revisit the spectral decomposition of the lamplighter random walk on the discrete circle, describing the spectral decomposition in terms of irreducible representations of the group $C_2\wr C_n$. In a similar way, the lamplighter on the complete graph is revisited in Section 10, now using the action of the hyperoctahedral group $C_2\wr S_n$.\\

In our join paper with T. Ceccherini-Silberstein \cite{CST2}, we analyzed several constructions that lead to multiplicity free permutation representations of wreath products. On the contrary, the harmonic analysis of the lamplighter random walk with a transitive group action leads to an example of a space with multiplicities. Moreover, in our examples the operators are not in the center of the commutant of the permutation representation, and therefore their diagonalization with irreducible eigenspaces requires a suitable explicit orthogonal decomposition of each isotypic component (see Proposition \ref{multspaces}). An example that leads to an operator in the center is in \cite{Schoolfield}; see also Remark \ref{lastrem}.

\section{Preliminaries and notation.}
If $X$ is a finite set, $L(X)$ will denote the space of all complex functions defined on $X$.
The space $L(X)$ will be endowed with the scalar product $\langle f_1,f_2\rangle_{L(X)}=\sum_{x\in X}f_1(x)\overline{f_2(x)}$, $f_1,f_2\in L(X)$. The symbol $\delta_x$ will denote the Dirac function centered at $x\in X$ and if $A\subseteq X$ then $\mathbf{1}_A$ is the characteristic function of $A$. Let $Y$ be another set. We will use the isomorphism $L(X\times Y)\cong L(X)\otimes L(Y)$, where for $f_1\in L(X)$,$f_2\in L(Y)$, $x\in X$ and $y\in Y$, we have $(f_1\otimes f_2)(x,y)=f_1(x)f_2(y)$. Let $C_2=\{0,1\}$ be the two elements cyclic group written additively. Then the set $\{0,1\}^X$ will denote the finite abelian group of all functions $\theta:X\rightarrow \{0,1\}$, with addition $(\theta+\omega)(x)=\theta(x)+\omega(x)$ mod 2. The identity of this group will be denoted by $\mathbf{0}_X$ (that is $\mathbf{0}_X\equiv 0$ on all $X$).

For $\theta, \omega \in \{0,1\}^X$, define the scalar product $\theta\cdot \omega = \sum_{x \in X}\theta(x) \omega(x)$ and set $\chi_\theta(\omega) = (-1)^{\theta\cdot \omega}$.
Then $\chi_\theta$ is a character of $C_2$ and the dual group is $\widehat{C_2^X}=\{\chi_\theta:\theta\in C_2^X\}$. We also recall the orthogonality relations $\langle\chi_\theta,\chi_\omega\rangle_{L(\{0,1\}^X)}=\frac{1}{2^{\lvert X\rvert}}\delta_{\theta,\omega}$.

Let $(X,E)$ be a graph simple, unoriented and without loops. We will think of the edge set $E$ as a subset of $\{\{x,y\}:x,y\in X,x\neq y\}$ and we will write $x \sim y$ to denote that $\{x,y\}$ is an edge. By  $\deg(x) = |\{y \in X: x \sim y\}|$ we will denote the degree of $x\in X$. The Markov operator of the graph is the linear selfadjoint operator $M:L(X)\rightarrow L(X)$ defined by setting

\[
(Mf)(x) = \frac{1}{\deg(x)}\sum_{y\sim x}f(y),
\]

while the adjacency operator is given by

\[
(Af)(x) = \sum_{y\sim x}f(y),
\]

for any $f\in L(X)$. If $X$ is regular of degree $k$, we have $M=\frac{1}{k}A$, but if $X$ is not regular, in general $M$ and $A$ have a different spectral theory. For instance, for the path the adjacent spectrum requires a discrete sine transform \cite{Biggs, Sca}, while its Markov spectrum requires a discrete cosine transform \cite{BDPX, Feller}.

If $g_1,g_2,\dotsc, g_m$ belong to a group $G$, then $\langle g_1,g_2,\dotsc g_m\rangle$ will denote the subgroup generated by $g_1,g_2,\dotsc,g_m$; if $v_1,v_2,\dotsc v_m$ belong to a vector space $V$, then $\langle v_1,v_2\dotsc v_m\rangle$ will denote the subspace spanned by $v_1,v_2\dotsc v_m$. If $G$ is a finite group, $(\rho,V)$ a unitary representation of $G$ and

\begin{equation}\label{isot}
V=\oplus_{j\in J}m_jW_j
\end{equation}

\noindent
is the decomposition of $V$ into irreducible representations $W_j$ where $m_jW_j=W_j\oplus\dotsc\oplus W_j$ $m_j$-times and $W_i,W_j$ are inequivalent for $i\neq j$, then we say that the $\{m_jW_j:j\in J\}$ are the {\it isotypic components} of $V$.

\section{Vertex lamplighter random walks.}\label{s, lumpvertex}

Let $(X,E)$ be a finite graph. Set

\[
 \mathcal{L}(X) = \left\{(\omega,x): \omega\in \{0,1\}^X, x \in X\right\}\equiv \{0,1\}^X\times X.
\]

Following \cite{Pe-Re}, we define a graph structure on  $\mathcal{L}(X)$ by declaring
two vertices $(\omega,x)$, $(\theta,y) \in \mathcal{L}(X)$ adjacent if
$x\sim y$ (in $X$) and $\omega(z) = \theta(z)$ for all $z \neq x,y$.
In other words, $x$ must be connected  to  $y$ and $\omega$  must take the same values of $\theta$ on $X\setminus\{x,y\}$.
The vertex lamplighter process on $X$ is the simple random walk on $\mathcal{L}(X)$.
Note that $L(\mathcal{L}(X)) \equiv L(\{0,1\}^X) \otimes L(X)$ and that
the Markov operator on $\mathcal{L}(X)$ is:

\[
[\mathcal{M}_X(F\otimes f)](\omega,x) = \frac{1}{4 \deg(x)}\sum_{y \sim x}[F(\omega)+
F(\omega + \delta_x)+ F(\omega+ \delta_y)+F(\omega+\delta_x+ \delta_y)]f(y),
\]

where  $F \in L(\{0,1\}^X)$, $f \in L(X)$ and $(\omega,x) \in \mathcal{L}(X)$.

If we define $V_\theta = \{\chi_\theta\otimes f: f \in L(X)\}$,
then we have the orthogonal decomposition

\begin{equation}\label{orthdec}
L\left(\mathcal{L}(X)\right) = \bigoplus_{\theta \in \{0,1\}^X}V_\theta.
\end{equation}

Now we show how to reduce the spectral analysis of $\mathcal{M}$ to a series of eigenvalues problem on $X$.
For $\theta \in  \{0,1\}^X$, we set $X_\theta=\{x\in X:\theta(x)=0\}$. We define a linear operator $M_\theta: L(X) \to L(X)$ by setting, for $f \in L(X)$ and $x \in X$,
\[
(M_\theta f)(x) = \left\{\begin{array}{cl}\frac{1}{\deg(x)} \sum_{\substack{y \in X_\theta:\\ y \sim x}}f(y) &
\mbox{ if  $\theta(x) = 0$}\\
0 & \mbox{ if  $\theta(x) = 1$.}
\end{array}\right.\]

\begin{lemma}\label{mtheta}
If $f \in L(X)$ then
\[
\mathcal{M}_X(\chi_\theta \otimes f) = \chi_\theta \otimes M_\theta f.
\]
\end{lemma}

\begin{proof}
\[
\begin{split}
\left[\mathcal{M}_X(\chi_\theta \otimes f)\right](x,\omega) & = \frac{1}{4 \deg(x)}\sum_{y \sim x}[
\chi_\theta(\omega)+ \chi_\theta(\omega+ \delta_x)+\chi_\theta(\omega+\delta_y)+\chi_\theta(\omega+\delta_x+\delta_y)]f(y).
\end{split}
\]

As the term in squared brackets  equals
\[
\begin{split}
\chi_\theta(\omega)[1 + \chi_\theta(\delta_x)+ \chi_\theta(\delta_y)+ \chi_\theta(\delta_x+ \delta_y)] & =
\chi_\theta(\omega)[1 + (-1)^{\theta(x)}+ (-1)^{\theta(y)}+ (-1)^{\theta(x)+ \theta(y)}] = \\
& = \left\{\begin{array}{lc}
4 \chi_\theta(\omega) & \mbox{ if $\theta(x) = \theta(y) = 0$}\\
0 & \mbox{ otherwise,}
\end{array}\right.
\end{split}
\]

if $\theta(x) = 0$ we have that

\[
\left[\mathcal{M}_X(\chi_\theta \otimes f)\right](x,\omega) = \chi_\theta(\omega)\frac{1}{\deg(x)}
\sum_{\substack{y \in X_\theta:\\ y \sim x}}f(y),
\]
while if $\theta(x) = 1$ then $[\mathcal{M}_X(\chi_\theta\otimes f)](x,\omega) = 0$. Therefore,
$\mathcal{M}_X(\chi_\theta \otimes f) = \chi_\theta \otimes M_\theta f$.
\end{proof}

\begin{remark}
{\rm In other words, the lamplighter random walk is $C_2^X$-invariant (cf.\eqref{gensym}). Moreover, \eqref{orthdec} is the decomposition of $L(\mathcal{L}(X))$ into irreducible $C_2^X$ representations, that is $V_\theta$ is the isotypic component corresponding to the character $\chi_\theta$. Then each $V_\theta$ is $C_2^X$-invariant and Lemma \ref{mtheta} is just the expression of the restriction of $\mathcal{M}$ to $V_\theta$. Lemma \ref{mtheta} may be also seen as a finite generalization of Lemma 3.7 in \cite{Ba-Woe}.}
\end{remark}

We now give a closer look at the operator $M_\theta$.
Let $X_\theta$ be as before and set $E_\theta = \{\{x,y\}\in E: x,y \in X_\theta\}$.
Clearly $X_\theta$, with edge set $E_\theta$, is a subgraph of $X$.
If we denote by $A_\theta$ the adjacency operator of $X_\theta$, then we have

\[
(M_\theta f)(x) = \frac{1}{\\deg(x)}(A_\theta f)(x).
\]

In particular, if $\deg(x) = k$ (i.e it is constant on $X$), then
$(M_\theta f)(x) = \frac{1}{k}(A_\theta f)(x)$ and therefore one can recover the spectrum of
$\mathcal{L}(X)$ by analyzing the adjacency spectra of  all the subgraphs of $X$ that may be obtained erasing some vertices  of $X$.

Let

\begin{equation}\label{e; spectra}
M_\theta  = \lambda_{\theta,1}P_{\theta,1} + \lambda_{\theta,2}P_{\theta,2}+ \cdots +
\lambda_{\theta, h(\theta)}P_{\theta, h(\theta)}
\end{equation}

be the spectral decomposition of $M_\theta$. That is,
$\lambda_{\theta,1}, \lambda_{\theta,2}, \ldots, \lambda_{\theta,h(\theta)}$ are the distinct nonzero eigenvalues
and $P_{\theta,j}$ is the orthogonal projection
of $L(X)$ onto  the eigenspace of $\lambda_{\theta,j}$. Clearly, if $X_\theta \subsetneq X$, $M_\theta$ has also
the eigenspace $L(X\setminus X_\theta)$, with eigenvalue equal to zero; this is omitted in \eqref{e; spectra}.

Let $Q_\theta: L(\mathcal{L}(X)) \to V_\theta$ be the orthogonal projection onto $V_\theta$. Then,
for $F = F(\omega,x) \in  L(\mathcal{L}(X))$, we have

\[
Q_\theta F = \chi_\theta\otimes \widetilde{Q}_\theta F,
\]

where $(\widetilde{Q}_\theta F)(x) = \frac{1}{2^{|X|}}\sum_{\omega \in \{0,1\}^X}F(\omega,x)\chi_\theta(\omega)$.

\begin{lemma}
The spectral decomposition of the operator $\mathcal{M}_X$ is given by:

\[
\mathcal{M}_X = \sum_{\theta \in \{0,1\}^X}\sum_{j = 1}^{h(\theta)}\lambda_{\theta,j}\left(\chi_\theta \otimes P_{\theta,j}
\widetilde{Q}_\theta\right),
\]

where ($\chi_\theta \otimes P_{\theta,j}\widetilde{Q}_\theta)F = \chi_\theta \otimes P_{\theta,j}\widetilde{Q}_\theta F$
for any $F \in L(\mathcal{L}(X))$.
The zero eigenvalues (in particular those corresponding to the space $L(X\setminus X_\theta)$) are omitted and the eigenvalues
$\{\lambda_{\theta,j}: \theta\in \{0,1\}^X, j = 1,2,\ldots, h(\theta)\}$ are not necessarily distinct.
\end{lemma}

\begin{proof}
It is obvious: if $F \in L(\mathcal{L}(X))$ then

\[
\begin{split}
\mathcal{M}_XF &= \mathcal{M}_X\left(\sum_{\theta \in \{0,1\}^X}\chi_\theta\otimes \widetilde{Q}_\theta F\right)=\\
& = \sum_{\theta \in \{0,1\}^X}\chi_\theta\otimes M_\theta\widetilde{Q}_\theta F =\\
& =  \sum_{\theta \in \{0,1\}^X}\sum_{j = 1}^{h(\theta)}\lambda_{\theta,j}
\left(\chi_\theta\otimes P_{\theta,j}\widetilde{Q}_\theta F\right).
\end{split}
\]

\end{proof}

\begin{remark}{\rm In other words, if $V_{\theta,j}$ is the eigenspace of $M_\theta$ corresponding to  $\lambda_{\theta,j}$,
$j = 0,1, \ldots,h(\theta)$ (with $V_{\theta,0}$ the eigenspace of $\lambda_{\theta,0}= 0$) then
$W_{\theta,j} = \{\chi_\theta \otimes  f:f \in V_{\theta, j}\}$ is the eigenspace of $M_\theta$ corresponding to
$\lambda_{\theta,j}$.}
\end{remark}

\begin{corollary}[$k-$step iterate]\label{sst}
The probability of going  from $(\omega,x)$ to $(\eta,y)$ in $k$ steps is equal to
\begin{equation}\label{e;nstepit}
\sum_{\theta \in \{0,1\}^X}\sum_{j = 1}^{h(\theta)}(\lambda_{\theta,j})^k \cdot\frac{\chi_\theta(\omega)\cdot \chi_\theta(\eta)}{2^{|X|}}(P_{\theta,j}\delta_x)(y).
\end{equation}
\end{corollary}

\begin{proof}
We have
\[
\widetilde{Q}_\theta(\delta_\omega\otimes \delta_x) =
\frac{1}{2^{|X|}}\chi_\theta(\omega)\delta_x
\]
and therefore
\begin{equation}\label{e;nstepit2}
\mathcal{M}_X^k(\delta_\omega\otimes \delta_x) =
 \sum_{\theta\in \{0,1\}^X}\sum_{j = 1}^{h(\theta)}(\lambda_{\theta,j})^k \cdot \frac{\chi_\theta(\omega)}{2^{|X|}}\chi_\theta \otimes P_{\theta,j}\delta_x
\end{equation}
from which \eqref{e;nstepit} follows immediately, since it is equal to
$[\mathcal{M}_X^k(\delta_\omega\otimes \delta_x)](\eta,y)$.
\end{proof}

Now we give the lamplighter version of the celebrated upper bound lemma of Diaconis and Shahshahani \cite{Diaconis}.

\begin{corollary}[Upper bound lemma]
Suppose that $X$ is connected. Assuming that $P_{{\bf{0}}_X,1}$  is the orthogonal projector on the space of constant value
functions, we have

\[
\begin{split}
\left\lVert[\mathcal{M}_X^k(\delta_\omega\otimes \delta_x) -\frac{1}{2^{|X|}|X|}
{\bf{1_{\mathcal{L}(X)}}} \right\rVert^2_{TV} & \leq
|X|\left\{\sum_{\substack{\theta\in \{0,1\}^X:\\ \theta \neq {\bf{0}}_X}}\sum_{j=1}^{h(\theta)}
|\lambda_{\theta,j}|^{2k}\|P_{\theta,j}\delta_x\|_{L(X)}^2 + \right.\\
&+ \left.\sum_{j=2}^{h({\bf{0}}_X)}|\lambda_{\mathbf{0}_X,j}|^{2k}\|P_{{\bf{0}}_X,j}\delta_x\|_{L(X)}^2\right\}.
\end{split}
\]

\end{corollary}

\begin{proof}

From \eqref{e;nstepit2} and the orthogonality relations for the characters $\chi_\theta$'s we get:

\[
\begin{split}
\left\lVert[\mathcal{M}_X^k(\delta_\omega\otimes \delta_x)-\frac{1}{2^{|X|}|X|}
{\bf{1_{\mathcal{L}(X)}}}\right\rVert^2_{L(\mathcal{L}(X))} & =
\sum_{\substack{\theta\in \{0,1\}^X:\\ \theta \neq {\bf{0}}_X}}
\sum_{j=1}^{h(\theta)}\frac{|\lambda_{\theta,j}|^{2k}}{2^{2|X|}}\|\chi_\theta \otimes P_{\theta,j}\delta_x\|_{L(\mathcal{L}(X))}^2 \\
&+\sum_{j=2}^{h(\mathbf{0}_X)}
\frac{|\lambda_{\mathbf{0}_X,j}|^{2k}}{2^{2\lvert X\rvert}}
\|\chi_{\mathbf{0}_X} \otimes P_{\mathbf{0}_X,j}\delta_x\|_{L(\mathcal{L}(X))}^2
=\\
& = \sum_{\substack{\theta\in \{0,1\}^X:\\ \theta \neq {\bf{0}}_X}}\sum_{j=1}^{h(\theta)}\frac{|\lambda_{\theta,j}|^{2k}}{2^{|X|}} \|P_{\theta,j}\delta_x\|_{L(X)}^2+\\
&+\sum_{j=2}^{h(\mathbf{0}_X)}
\frac{|\lambda_{\mathbf{0}_X,j}|^{2k}}{2^{\lvert X\rvert}}
\|P_{\mathbf{0}_X,j}\delta_x\|_{L(X)}^2.
\end{split}
\]

Then the upper bound lemma follows immediately from the Cauchy-Schwarz inequality.
\end{proof}

\begin{remark}{\rm The hypothesis that $X$ is connected guarantees that the multiplicity
 of the eigenvalue $1$ is equal to 1.}
\end{remark}

\section{The vertex lamplighter random walk on the complete graph. }\label{completegraph}

Suppose that $(X,E)$ is the complete graph on $n$ vertices. We identify $X$ with $\{1,2,\dotsc,n\}$. Now for any $\theta\in\{0,1\}^X$, the graph $(X_\theta,E_\theta)$ is the complete graph on $\lvert X_\theta\rvert$ vertices. We recall that the eigenspaces of the adjacency operator on the complete graph on $m$ vertices are the space of constant functions and its orthogonal complement, with corresponding eigenvalues $m-1$ and $-1$.\\

For any $\theta \in \{0,1\}^X$, define the projector $P_\theta:L(X)\rightarrow L(X)$ by setting

\[
P_\theta f(x) = \left\{\begin{array}{cl}
\frac{1}{\lvert X_\theta \rvert}\sum_{y \in X_\theta}f(y) & \mbox{if $x \in X_\theta$}\\
0 & \mbox{if $x \not\in X_\theta$,}
\end{array}\right. \qquad \text{\rm for any}\quad f\in L(X).
\]

For $\lvert X_\theta\rvert>1 $, the spectral decomposition of the operator $M_\theta$  is given by

\[
M_\theta = \frac{\lvert X_\theta\rvert-1}{n-1}P_\theta - \frac{1}{n-1}(R_\theta-P_\theta)
\]

where $R_\theta:L(X)\rightarrow L(X_\theta)$ is the orthogonal
projection from $L(X)$ onto $L(X_\theta)$. In the notation
introduced above, we have: $h(\theta) = 2$, $\lambda_{\theta,1} =
\frac{\lvert X_\theta\rvert -1}{n-1}$, $\lambda_{\theta,2} =
-\frac{1}{n-1}$, $P_{\theta,1} = P_\theta$ and $P_{\theta,2} =
R_\theta-P_\theta$.

Clearly, if $\lvert X_\theta\rvert = 1$, then $X_\theta = \{x\}$ for some $x \in X$ and $M_\theta \equiv 0$.

If $x \in X_\theta$ we have
\[
(P_\theta \delta_x)(y) = \left\{\begin{array}{ll}
\frac{1}{\lvert X_\theta\rvert} & \mbox{if $y \in X_\theta$}\\
0 & \mbox{if $y \notin X_\theta$}
\end{array}
\right.
\]
and therefore, $\|P_\theta \delta_x\|^2_{L(X)} = \frac{1}{\lvert X_\theta\rvert}$ and
$\|(R_\theta-P_\theta) \delta_x\|^2_{L(X)} = \frac{\lvert X_\theta\rvert-1}{\lvert X_\theta\rvert}$.

If $x \notin X_\theta$ then $\mathcal{M}_X(\delta_\omega \otimes  \delta_x) \equiv 0$.

Denote by
\[
\mathcal{A}_i = \{\theta \in \{0,1\}^X: \lvert X_\theta\rvert = i+1\}
\]
and observe that $|\mathcal{A}_i| = \binom{n}{i+1}$.
Now we are in position to estimate the rate of convergence to the stationary distribution.

\begin{proposition}
There exists $C>0$ such that if   $k \geq \frac{n}{2}(\log n+c)$ with $c\geq 0$, we have
\[
\left\|\mathcal{M}_X^k(\delta_{\omega_0} \otimes \delta_y) -\frac{1}{2^n n}
{\bf{1_{\mathcal{L}(X)}}}\right\|^2_{TV} \leq C\exp({-c}).
\]
\end{proposition}
\begin{proof}

By the upper bound lemma, we have

\begin{multline}\label{cs}
\left\|\mathcal{M}_X^k(\delta_{\omega_0} \otimes \delta_{x_0}) -\frac{1}{2^n n}
{\bf{1_{\mathcal{L}(X)}}}\right\|^2_{TV} \leq\\
\leq |X|\left\{\sum_{i = 1}^{n-2}\sum_{\theta\in \mathcal{A}_i}\sum_{j = 1}^{2}
|\lambda_{\theta,j}|^{2k}\|P_{\theta,j}\delta_{x_0}\|_{L(X)}^2 +
|\lambda_{\mathbf{0}_X,2}|^{2k}\|P_{\mathbf{0}_X,2}\delta_{x_0}\|_{L(X)}^2\right\}\leq\\
\leq n\left\{
 \sum_{i = 1}^{n-2}\left[\binom{n}{i+1}\left(\frac{i}{n-1}\right)^{2k}\frac{1}{i+1}\right]\right.\\
 +\left.\sum_{i = 1}^{n-2}\left[\binom{n}{i+1}\left(\frac{1}{n-1}\right)^{2k}\frac{i}{i+1}\right] +\left(\frac{1}{n-1}\right)^{2k}\frac{n-1}{n}
\right\}.
\end{multline}

Note that in the third step we have an inequality as we have to take into account the cases when $x_0 \notin X_\theta$.\\

The largest nontrivial eigenvalue is $\frac{n-2}{n-1}$ and the
corresponding term in \eqref{cs} is
$\frac{n^2}{n-1}\left(1-\frac{1}{n-1}\right)^{2k}<\frac{n^2}{n-1}\exp\left(-\frac{2k}{n-1}\right)$,
which becomes $<1$ when
$k>\frac{n-1}{2}\log\frac{n^2}{n-1}\sim\frac{n}{2}\log n$. It
remains to show that the other part of \eqref{cs} goes to zero
faster.

Suppose that $k = \frac{1}{2}n(\log n+c)$ with $c>0$.
The last term in (\ref{cs}) is clearly smaller than $e^{-c}$ if $n$ is sufficiently large.
Moreover, it is obvious that the second sum is dominated by the first sum, and therefore we are left to estimate the first sum. With the change of variable $i\rightarrow n-i-1$, we have:

\[
\begin{split}
\sum_{i = 1}^{n-2}\binom{n}{n-1- i}\left(\frac{i}{n-1}\right)^{2k}\frac{n}{i+1} &
= \sum_{i = 1}^{n-2}\binom{n}{i}\left(1-\frac{i}{n-1}\right)^{2k}\frac{n}{n-i} \leq\\
& \leq \sum_{i = 1}^{n-2}\frac{n^i}{i!}\exp\left(-\frac{2ki}{n-1}\right)\frac{n}{n-i}\leq \\
\text{setting}\quad k=\frac{n}{2}(\log n+c)\qquad& \leq \sum_{i = 1}^{n-2}\exp\left(-\log(i!)-ic+\log\frac{n}{n-i}\right) \leq \\
& \leq\exp(-c)\sum_{i = 1}^{n-2}\exp\left[-i\log(i)+i-1+\log n-\log(n-i)\right],
\end{split}
\]

since $\log(i!)\geq i\log i-i+1$ and $-ic\leq -c$.
In order to complete the proof, we just need to bound the last sum by a constant independent of $n$.
Observe that $i \mapsto h(i) =[-i\log(i)+i-1+\log n-\log(n-i)]+i$ has derivative equal to $-\log(i)+1 +\frac{1}{n-i}$, which is negative if $i\geq 10$. Moreover, if $n\geq 20$ then $h(10)\leq 0$ and therefore we can conclude that

\begin{multline}
\sum_{i = 1}^{n-2}\exp\left(-i\log(i)+i-1+\log n-\log(n-i)\right) \\
\leq \sum_{i = 1}^{10}\exp\left(-i\log(i)+i-1+\log n-\log(n-i)\right) +
\sum_{i = 11}^{+\infty}\exp(-i) \leq C.
\end{multline}
 \end{proof}

Now we give the corresponding lower bound, showing that the random walk has a cut-off at $k=\frac{n}{2}\log n$.

\begin{proposition}
Let $p^{(k)} = \mathcal{M}_X^k(\delta_{\omega_0} \otimes \delta_{x_0})$ be the probability after $k$ steps  starting from the point $(\omega_0,x_0)$ and let $\pi$ be the uniform distribution. Then
for $k = \frac{1}{2}n(\log n-c)$, $0<c<\log n$ and $n$ large we have
\[
\|p^{(k)}-\pi\|_{TV} \geq 1 - 20 e^{-c}.
\]
\end{proposition}
\begin{proof}
For $x = 1,2,\ldots, n$ let $f_x:X\to \mathbb{C}$ be the characteristic function  of $X_{\delta_x}$. Moreover, for $x\neq y$  set $\theta_{x,y}= \delta_x+ \delta_y$ and let
$f_{x,y}$ be the characteristic function of $X_{\theta_{x,y}}$.
We have the following equality:

\begin{equation}\label{lampiolb}
\left(\chi_{\delta_x}\otimes  f_x\right)\left(\chi_{\delta_y}\otimes  f_y\right)=
\left\{\begin{array}{ll}
\chi_{{\bf{0}}_X}\otimes f_x & \mbox{ if } x = y\\
\chi_{\theta_{x,y}} \otimes f_{x,y}& \mbox{ if } x \neq  y.
\end{array}
\right.
\end{equation}

Given a probability distribution $p$ on $\mathcal{L}(X)$ and a function $f:\mathcal{L}(X)\to \mathbb{C}$, the expected value of $f$ with respect to $p$ is
$E_p(f) = \sum_{(\omega,x) \in \mathcal{L}(X)}p(\omega,x)f(\omega,x)$,
while the variance of $f$ is
$Var_p(f) = E_p(f^2) - E_p(f)^2$.

Consider the function
\[
F= \chi_{\delta_1}\otimes f_1+ \chi_{\delta_2}\otimes f_2+ \cdots + \chi_{\delta_n}\otimes f_n
\]

which is an eigenvector of the operator $\mathcal{M}_X$, with eigenvalues $\frac{n-2}{n-1}$.
In what follows, we suppose that $\omega_0=\mathbf{0}_X$; this implies
$F(\omega_0,x_0)=n-1$.
In virtue of (\ref{lampiolb}), we have

\begin{equation}\label{F2}
F^2 = \chi_{{\bf{0}}_X}\otimes f_1+ \cdots + \chi_{{\bf{0}}_X}\otimes f_n+ \sum_{x \neq y}\chi_{\theta_{x,y}}\otimes f_{x,y} = (n-1)\chi_{{\bf{0}_X}}\otimes {\bf{1}}_X + \sum_{x \neq y}\chi_{\theta_{x,y}}\otimes f_{x,y}.
\end{equation}

But $E_{p^{(k)}}(f)=[\mathcal{M}^k_X(f)](\omega_0,x_0)$ and therefore

\begin{equation}\label{e;cumemu}
E_{p^{(k)}}(F) = \left(\frac{n-2}{n-1}\right)^k\left(f_1(x_0)+ f_2(x_0) + \cdots +f_n(x_0)\right) =  (n-1)\left(\frac{n-2}{n-1}\right)^k.
\end{equation}

Similarly, by (\ref{F2}) we have,
\[
E_{p^{(k)}}(F^2)
=  n-1 + \left(\frac{n-3}{n-1}\right)^k\sum_{x\neq y}f_{x,y}(x_0) = n-1 + (n-1)(n-2)\left(\frac{n-3}{n-1}\right)^k.
\]

and therefore

\begin{equation}\label{e;cuvamu}
Var_{p^{(k)}}(F) =  n-1 + {(n-1)(n-2)}\left(\frac{n-3}{n-1}\right)^k-  (n-1)^2\left(\frac{n-2}{n-1}\right)^{2k}\leq n-1.
\end{equation}

Since $\pi = \frac{1}{2^n n}{\bf{1_{\mathcal{L}(X)}}}$ is the uniform distribution, we have $E_\pi(F) = 0$ and $Var_\pi(F) = n-1$.

Now define ${\bf{A}}_\beta = \{(\omega,x) \in \mathcal{L}(X): |F(\omega,x)|<\beta\sqrt{n-1}\}$, where $\beta$ is a constant $0< \beta <
\frac{1}{\sqrt{n-1}}E_{p^{(k)}}(F)$ that  will be suitably chosen later.
From Markov's  inequality  it follows that

\begin{equation}\label{e;cuA}
\begin{split}
\pi({\bf{A}}_\beta ) & =1 -\pi\{(\omega,x):|F(\omega,x)| \geq \beta\sqrt{n-1}\}\geq\\
& \geq 1 - \frac{1}{\beta^2(n-1)}E_\pi(F^2) = 1 - \frac{1}{\beta^2}.
\end{split}
\end{equation}

In the same way,  from Chebyshev's inequality
and the fact that ${\bf{A}}_\beta \subseteq \{(\omega,x) \in \mathcal{L}(X): |F(\omega,x)-E_{p^{(k)}}(F)|\geq  E_{p^{(k)}}(F) - \beta \sqrt{n-1}\}$, we have

\begin{equation}\label{e;cufe}
p^{(k)}({\bf{A}}_\beta) \leq
 \frac{Var_{p^{(k)}}(F)}{(E_{p^{(k)}}(F)-\beta\sqrt{n-1})^2}.
\end{equation}

 Set $k=\frac{n}{2}(\log n -c)$, $0<c<\log n$. From the Taylor expansion of the logarithm, it follows that $\log(1-t)=-t-\frac{t^2}{2}\eta(t)$,
with $\eta(t)\geq 0$ and $\lim_{t\rightarrow 0}\eta(t)=1$. Applying this asymptotic expansion
 to the right hand side of (\ref{e;cumemu}), we get

\begin{equation*}
\begin{split}
E_{p^{(k)}}(F)
=&(n-1)\exp\left\{\left[-\frac{1}{n-1}-\frac{1}{2(n-1)^2}\cdot\eta\left(\frac{1}{n-1}\right)\right]\cdot\frac{n}{2}(\log n-c) \right\} = \\
=&\frac{n-1}{\sqrt{n}}e^{c/2}
\exp\left\{\frac{c-\log n}{2(n-1)}\left[1+\frac{n}{2(n-1)}\cdot\eta\left(\frac{1}{n-1}\right)\right]\right\}
\end{split}
\end{equation*}

and therefore for $n$ large we have

\begin{equation}\label{e;cuscf}
E_{p^{(k)}}(F)\geq  \frac{3}{4}\sqrt{n-1}e^{c/2}
\end{equation}

Choosing $\beta =  \frac{e^{c/2}}{2}$ and taking in account (\ref{e;cuvamu}) and (\ref{e;cuscf}),
we have that (\ref{e;cufe}) becomes
\begin{equation}\label{e;cuB}
p^{(k)}({\bf{A}}_\beta) \leq \frac{n-1}{(\frac{3}{2}\beta\sqrt{n-1}-\beta\sqrt{n-1})^2} = \frac{4}{\beta^2}.
\end{equation}

and therefore

\[
\lVert p^{(k)}-\pi\rVert_{TV}\geq\pi({\bf{A}}_\beta)-p^{(k)}({\bf{A}}_\beta)\geq 1-\frac{5}{\beta^2}=1-20e^{-c}.
\]

\end{proof}

\section{Edge lamplighter random walks.}
Let $(X,E)$ be again a finite graph.
Set $\mathcal{L}(E) = \{(\omega,x): \omega \in \{0,1\}^E, x \in X\}$ and  define a graph structure on
$\mathcal{L}(E)$ by declaring two vertices $(\omega, x)$, $(\theta,y) \in \mathcal{L}(E)$ adjacent when
$x \sim y$ and $\omega(e) = \theta(e)$, for all $e \in E\setminus \{\{x,y\}\}$. Therefore $x$ must be connected to $y$ and
$\omega$ must take the same values of $\theta$ on  any edge different from $\{x,y\}$.
The simple random walk on $\mathcal{L}(E)$ is the following: the lamplighter moves from a vertex $x$
to an adjacent vertex $y$ with equal probability; when he moves from the vertex $x$ to the vertex $y$ he changes randomly
 the state of the lamp on the edge $\{x,y\}$.

The Markov operator on $\mathcal{L}(E)$ is

\[
[\mathcal{M}_E(F\otimes f)](\omega, x)= \frac{1}{2 \deg(x)}\sum_{\substack{y \in X:\\y \sim x}}\left[F(\omega)+ F(\omega+ \delta_{\{x,y\}})\right]f(y)
\]

for $F \in L(\{0,1\}^E)$, $f \in L(X)$, $x \in X$  and
$\omega \in \{0,1\}^E$.
Clearly $L(\mathcal{L}(E)) = L(\{0,1\}^E)\otimes L(X)$.

For $\theta, \omega \in \{0,1\}^E$,
we define $V_\theta = \{\chi_\theta\otimes f: f \in L(X)\}$, where $\chi_\theta$ is again the character associated to $\theta$; we have the orthogonal decomposition

\[
L\left(\mathcal{L}(E)\right) = \bigoplus_{\theta \in \{0,1\}^E}V_\theta.
\]

For $\theta \in  \{0,1\}^E$, define  the linear operator $M_\theta: L(X) \to L(X)$ by setting

\[
(M_\theta f)(x) =
\frac{1}{\deg(x)} \sum_{\substack{y \in X:\\ y \sim x\\ \theta(\{x,y\}) = 0}}f(y).
\]

\begin{lemma}
If $f \in L(X)$ then
\[
\mathcal{M}_E(\chi_\theta \otimes f) = \chi_\theta \otimes M_\theta f.
\]
\end{lemma}
\begin{proof}
For $x \in X$ and $\omega \in \{0,1\}^E$, we have

\[
\begin{split}
[\mathcal{M}_E(\chi_\theta \otimes f)](x,\omega)  & = \frac{1}{2 \deg(x)}
\sum_{\substack{y \in X:\\ y \sim x}}\left[\chi_\theta(\omega)+ \chi_\theta(\omega)
 (-1)^{\theta(\{x,y\})}\right]f(y) = \\
& = \chi_\theta(\omega) \frac{1}{2 \deg(x)} \sum_{\substack{y \in X:\\ y \sim x}}2(1 - \theta(\{x,y\})f(y)= \\
& = \chi_\theta(\omega)(M_\theta f) (x).
\end{split}
\]

The second step follows from   the observation that $1 + (-1)^\epsilon = 2(1-\epsilon)$ if $\epsilon\in\{0,1\}$.

\end{proof}

We now analyze  the operator $M_\theta$ more closely.
Set $E_\theta = \{ e \in E: \theta(e) = 0\}$. Then

\[
(M_\theta f)(x) = \frac{1}{\deg(x)} (A_\theta f)(x),
\]

where $A_\theta$ is the adjacency operator of the graph $(X, E_\theta)$. Note that $(X,E_\theta)$ is obtained from $(X,E)$ by deleting the edges $\{x,y\}$ such that $\theta(\{x,y\})=1$.
In particular, if $X$ is regular,  $\deg(x) = k$ and
$M_\theta = \frac{1}{k}A_\theta$.
As in Section \ref{s, lumpvertex}, let
$M_\theta  = \lambda_{\theta,1}P_{\theta,1} + \lambda_{\theta,2}P_{\theta,2}+ \cdots +
 \lambda_{\theta, h(\theta)}P_{\theta, h(\theta)}$
be the spectral decomposition of the operator $M_\theta$. Arguing as in Section \ref{s, lumpvertex}, one can get the spectral decomposition
of $\mathcal{M}_X$ in the form

\[\mathcal{M}_X = \sum_{\theta \in \{0,1\}^E}\sum_{j = 1}^{h(\theta)}
\lambda_{\theta,j}(\chi_\theta\otimes P_{\theta,j}\widetilde{Q}_\theta),
\]

with $W_{\theta,j} =
 \{\chi_\theta\otimes f : f \in V_{\theta,j}\}$ the eigenspace
corresponding to $\lambda_{\theta,j}$. In particular, now we have:

\begin{proposition}[Upper bound lemma II]
Suppose that $X$ is connected. Assuming that $P_{{\bf{0}}_X,1}$  is the orthogonal projector on the space of constant value
functions, we have
\[
\begin{split}
\left\|[\mathcal{M}_X^k(\delta_\omega\otimes \delta_x) -\frac{1}{2^{|E|}|X|}
{\bf{1_{\mathcal{L}(E)}}} \right\|^2_{TV} & \leq
\lvert X\rvert \left\{\sum_{\substack{\theta\in \{0,1\}^E:\\ \theta \neq {\bf{0}}_X}}\sum_{j=1}^{h(\theta)}
|\lambda_{\theta,j}|^{2k}\|P_{\theta,j}\delta_x\|_{L(X)}^2\right. + \\
&\left.+\sum_{j=2}^{h({\bf{0}}_X)}|\lambda_{0,j}|^{2k}\|P_{{\bf{0}}_X,j}\delta_x\|_{L(X)}^2\right\}.
\end{split}
\]
\end{proposition}

\begin{remark}
{\rm In general, the explicit diagonalization of all the operators
$M_\theta$ is quite a difficult (or impossible) task. For instance,
if $X$ is the complete  graph on $n$ vertices, it requires the
knowledge of the adjacency  spectrum of all graphs on $k \leq n $
vertices. Examples of graphs for which this is feasible are: the
path, the star and the discrete circle. In the following section, we
analyzed the edge lamplighter random walk on the discrete circle.
The path is analyzed in \cite{Sc-To}, also using Radon transforms on
a finite trees and a finite analogous of the construction of
Bartholdi and Woess \cite{Ba-Woe}.}
\end{remark}

\section{The edge lamplighter random walk on the discrete circle.}\label{disccircl}

Let $C_n$ be the discrete circle on $n$ points, that is the graph
with vertex set $C_n=\{0,1,2,\cdots,n-1\}$ and edge set
$E_n=\{\{0,1\},\{1,2\},\dotsc,\{n-2,n-1\},\{n-1,1\}\}$. Let $P_n$ be
the path of length $n-1$, that is the graph with vertex set
$P_n=\{0,1,2,\cdots,n-1\}$ and edge set
$\{\{0,1\},\{1,2\},\dotsc,\{n-2,n-1\}\}$. Consider the lamplighter
random walk on $C_n$, with the lamps on the edges. Clearly, if we
delete some edges of $C_n$, the resulting graph consists of a series
of disjoint paths; if we do not delete any edge, then we are
considering $C_n$ itself. We need two elementary facts of discrete
Fourier analysis; see \cite{M-P,Strang} for more details.

Let $A_n$ be the $n\times n$ circulant matrix

\[
A_n=\frac{1}{2}
\begin{pmatrix}
0       & 1       &0        &  \dotso      &             0&     1 \\
1       & 0       & 1       &              & &      0 \\
     \vdots   & \ddots & \ddots & \ddots       &      &\vdots \\
       0 &       &     & 1 & 0            & 1     \\
1 &0      & \dotso     & 0      & 1            & 0     \\
\end{pmatrix}.
\]

Set $w=\exp\left(\frac{2\pi i}{n}\right)$. Then the $n\times n$ symmetric matrix
\[
F_n = \frac{1}{\sqrt{n}}
 \begin{pmatrix}
1      & 1                  &1                & \dotso & 1              \\
1      & w^{-1}        & w^{-2}     & \dotso & w^{-(n-1)}  \\
1      & w^{-2}        & w^{-4}     & \dotso & w^{-2(n-1)}\\
\vdots & \vdots & \vdots    & \dotso          & \vdots    \\
1      & w^{-(n-1)}     & w^{-2(n-1)}  & \dotso &w^{-(n-1)(n-1)}\\
\end{pmatrix}
\]

is unitary and diagonalizes $A_n$:

\begin{equation}\label{circlespec}
F_nA_n\overline{F}_n=\begin{pmatrix}
\cos\frac{2\pi}{n} &&&\\
&\cos\frac{4\pi}{n}\\
&& \ddots&\\
&&& \cos \frac{2(n-1)\pi}{n}
\end{pmatrix}.
\end{equation}

Analogously, let $B_n$ be the $n\times n$ tridiagonal matrix

\[
B_n=\frac{1}{2}
\begin{pmatrix}
0       & 1       &         &        &              &       \\
1       & 0       & 1       &              &       \\
        & \ddots & \ddots & \ddots       &        \\
        &        &      1 & 0            & 1     \\
        &        &        & 1            & 0     \\
\end{pmatrix}.
\]

Then the $n\times n$ symmetric matrix

\[
S_n=\sqrt{\frac{2}{n+1}}
\begin{pmatrix}
\sin\frac{\pi}{n+1} & \sin \frac{2\pi}{n+1} &\dotso &\sin \frac{n\pi}{n+1}\\
\sin\frac{2\pi}{n+1} & \sin \frac{4\pi}{n+1} &\dotso  &\sin \frac{2n\pi}{n+1}\\
\vdots & \vdots &  & \vdots \\
\sin \frac{n\pi}{n+1} & \sin \frac{2n\pi}{n+1}&\dotso &\sin \frac{n^2\pi}{n+1}\\
\end{pmatrix}
\]

\noindent
is orthogonal and diagonalizes $B_n$:

\begin{equation}\label{specpath}
S_n B_n S_n=
\begin{pmatrix}
\cos\frac{\pi}{n+1} &&&\\
& \cos\frac{2\pi}{n+1} &&\\
&& \ddots&\\
&&& \cos \frac{n\pi}{n+1}
\end{pmatrix}.
\end{equation}

Clearly, \eqref{circlespec} is just the computation of the Markov spectrum of the circle, while \eqref{specpath} is just the computation of the $\frac{1}{2}$adjacency spectrum of the path $P_n$. In what follows, to simplify terminology, we will refer to \eqref{circlespec} and to \eqref{specpath} respectively as the spectrum of the circle and the spectrum of the path (note that, with this terminology, the spectrum of $P_2$ is $\{\pm\frac{1}{2}\}$).
The following theorem must be compared with the results of spectral analysis on finite trees in \cite{Strang2,Sca,Sc-To2}.

\begin{theorem}

The spectrum of the edge lamplighter random walk on $C_n$ is given
by:

\[
\{1\}\cup\{0\}\cup\left\{\cos\frac{h\pi}{k}:3\leq k\leq n+1, 1\leq h\leq k-1\;\text{\rm and}\;(h,k)=1\right\}.
\]

when $n$ is odd, and

\[
\{1\}\cup\{0\}\cup\{-1\}\cup\left\{\cos\frac{h\pi}{k}:3\leq k\leq n+1, 1\leq h\leq k-1\;\text{\rm and}\;(h,k)=1\right\}.
\]

when $n$ is even.

Moreover the multiplicities of the eigenvalues are the following.

\begin{enumerate}

\item Suppose that
$3\leq k\leq n+1$, $1\leq h\leq k-1$, $(h,k)=1$ and $n+1=kq+r$ with
$0\leq r \leq k-1$. Then the multiplicity of the eigenvalue $\cos\frac{h\pi}{k}$ is equal to

\begin{itemize}
\item $n\frac{2^n-2^{r-1}}{2^k-1}$ when $r\neq 0,1$;
\item $n\frac{2^n-1}{2^k-1}$ when $r= 1$ and $hq$ is odd;
\item $n\frac{2^n-1}{2^k-1}+2$ when $r= 1$ and $hq$ is even;
\item $n\frac{2^n+2^{k-1}-1}{2^k-1}$ when $r= 0$
\end{itemize}

\item The multiplicity of $1$ is always equal to 1; the multiplicity of $-1$ is equal to 1 when $n$ is even, and is equal to 0 when $n$ is odd.

\item The multiplicity of 0 is equal to

\begin{itemize}
\item $\frac{n}{3}2^n+\frac{n}{3}$ if $n$ is odd;
\item $\frac{n}{3}2^n-\frac{n}{3}$ if $n\equiv 2$ mod 4;
\item $\frac{n}{3}2^n-\frac{n}{3}+2$ if $n\equiv 0$ mod 4;
\end{itemize}

\end{enumerate}

\end{theorem}

\begin{proof}
We will say that $\theta\in\{0,1\}^{E_n}$ has a {\em segment of
length} $1\leq l\leq n-1$ if there exists $t\in C_n$ such that
$\theta(\{t-1,t\})=1$, $\theta(\{t,t+1\})=\theta(\{t+1,t+2\})=\dotso
=\theta(\{t+l-1,t+l\})=0$, $\theta(\{t+l,t+l+1\})=1$, where the
numbers $t, t+1,\dotsc,t+l$ are considered mod $n$; we will also say
that the segment is in position $t$. Clearly, there exist exactly
$2^{n-l-2}$ distinct $\theta$'s with a segment of length $l$ in
position $t$ (if $l=n-1$ there exists only one $\theta$), and
therefore any eigenvalue of $P_{l+1}$ appears $n2^{n-l-2}$ times as
an eigenvalue of the lamplighter random walk ($n$ times for
$l=n-1$). The problem is that the same number may be an eigenvalue
of $P_{l+1}$ for different values of $l$
and that it may be also an eigenvalue of $C_n$ (that corresponds to the case $l=n$.)\\

Consider the eigenvalue $\cos\frac{h\pi}{k}$, with $3\leq k\leq
n+1$, $1\leq h\leq k-1$, $(h,k)=1$. Suppose that $n+1=kq+r$ with
$2\leq r \leq k-1$. From \eqref{specpath} we deduce that
$\cos\frac{h\pi}{k}$ is an eigenvalue of any segment of length
$sk-2$, for $s=1,2,\dotsc,q$. Moreover, we cannot have $sk-2=n-1$
(because $r\neq 0$) and the eigenvalue cannot appear in the spectrum
of $C_n$ (because $r\neq 1$). Then the multiplicity is equal to:

\begin{equation}\label{molteigen}
\sum_{s=1}^q n2^{n-ks}=n2^n\sum_{s=1}^q\left(\frac{1}{2^k}\right)^s=n\frac{2^n-2^{r-1}}{2^k-1}.
\end{equation}

Now suppose that $r=1$, that is $n=qk$. Then we must consider also the spectrum of $C_n$. But we have $\cos\frac{h\pi}{k}=\cos\frac{2\pi j}{n}$, with $0<j<\frac{n}{2}$ if and only if $qh=2j$, that is $\cos\frac{h\pi}{k}$ appears as an eigenvalue of $C_n$ if and only if $qh$ is even. Moreover, $\cos\frac{2\pi j}{n}=\cos\frac{2\pi (n-j)}{n}$, and therefore any eigenvalue of $C_n$ different from $\pm 1$ has multiplicity two. Arguing as in \eqref{molteigen}, we immediately get the formulas for $r=1$ in the statement.\\

If $r=0$, that is $n+1=kq$, then we have just to correct \eqref{molteigen} (to consider the eigenvalue coming from the segments of length $n-1$): now it becomes $\sum_{s=1}^{q-1} n2^{n-ks}+n=n\frac{2^n+2^{k-1}-1}{2^k-1}$.\\

Clearly, $\pm 1$ are not ($\frac{1}{2}$ adjacency) eigenvalues of any segment; $1$ is always a multiplicity one eigenvalue of $C_n$ and $-1$ is a (multiplicity one) eigenvalue of $C_n$ if and only if $n$ is even.

It remains to prove the formulas for the multiplicity of the null
eigenvalue. First of all, note that $\cos\frac{j\pi}{k}=0$ exactly
when $k=2j$. Then any segment of even length yields a null
eigenvalue; if $n+1=2q+r$, with $0\leq r \leq 1$, then from the
segments we find

\begin{equation}
\begin{split}
\sum_{j=2}^q n2^{n-2j}=n\frac{2^n-4}{12}\qquad\text{if}\quad r=1\\
\sum_{j=2}^{q-1} n2^{n-2j}+n=n\frac{2^n+4}{12}\qquad\text{if}\quad r=0
\end{split}
\end{equation}

times the null eigenvalue. But the null eigenvalue arises also from the complements of the segments; that is, if $\theta(t)=0$, $\theta(t+1)=\theta(t+2)=\dotso =\theta(t+l)=1$, $\theta(t+l+1)=0$ then this part of $\theta$ yields $l-1$ times the null eigenvalue. Arguing as in \eqref{molteigen}, this way we get a total amount of

\[
\sum_{l=2}^{n-2}n2^{n-l-2}(l-1)+\underset{l=n-1}{n(n-2)}+\underset{l=n}{n}=\frac{n}{4}2^n
\]

times the null eigenvalue.
Finally, 0 is an eigenvalue of $C_n$ (with multiplicity 2) if and only if $n\equiv 0$ mod 4.

\end{proof}

\begin{remark}{\rm We recall that the $L^2$ (or chi square) distance between the distribution after $k$-steps and the stationary (in this case the uniform) distribution is just \cite{BDPX}

\[
2^{\vert E\rvert}\lvert X\rvert\left\lVert \mathcal{M}_E^k(\delta_\omega\otimes \delta_x)-\frac{1}{2^{\vert E\rvert}\lvert X\rvert} \mathbf{1}_{\mathcal{L}(E)}\right\rVert^2.
\]

For the lamplighter random walk on the discrete circle, the $L^2$ convergence to the stationary distribution is slower than the total variation convergence. This is shown in \cite{Pe-Re}, p.828. The first convergence requires order $n^3$ steps, while the second requires order $n^2$. A similar phenomenon is discussed in \cite{BDPX}. Using our spectral computations and the techniques in \cite{Diaconis}, it is easy to prove that the $L^2$ distance is bounded above by

\begin{equation}\label{upcircle}
2n\sum_{l=1}^{n-2}2^{n-l-2}\exp\left(-\frac{\pi^2k}{(l+2)^2}\right)+2n\exp\left(-\frac{\pi^2k}{(n+1)^2}\right)+\exp\left(-\frac{\pi^2k}{n^2}\right),
\end{equation}

which goes to zero exponentially after $k=\frac{n^2}{\pi^2}(n+c)$, $c>0$, steps. Note that in \eqref{upcircle} there is not a dominant term; the last term becoming  $<1$ is the term corresponding to $l+2=\frac{2}{3}n$, and this happens when $k=\frac{4\log 2}{27\pi^2}n^3$, but it is smaller than the term for $l=n-2$ when $k=\frac{\log 2}{2\pi^2}n^3$.}
\end{remark}

\section{A general decomposition for lamplighters on homogeneous spaces.}

In this section, we give a decomposition theorem in the case the graph is a homogeneous space. It is natural to prove this theorem in a more general form, that covers many other cases, such as the signed Bernoulli-Laplace diffusion model \cite{Schoolfield}.
Let $G$ be a finite group and $Z$ a finite homogeneous $G$-space. The group $G$ acts on $C_2^Z$ by setting, for $\omega\in C_2^Z$, $g\in G$ and $z\in Z$, $g\omega(z)=\omega(g^{-1}z)$.
The wreath product of $C_2$ by $G$ (with respect to the action of $G$ on $Z$)
is the set $C_2\wr G=\{(\omega,g):\omega\in C_2^Z, g\in G)\}\equiv C_2^Z\times G$ with the composition law: $(\theta,g)\cdot (\omega,h)=(\theta+g\omega,gh)$, for $\theta,\omega\in C_2^Z$, $g,h\in G$. The identity is given by: $(\mathbf{0}_Z,1_G)$, where $1_G$ is the identity of $G$; the inverse of an element is given by the formula: $(\theta,g)^{-1}=(g^{-1}\theta,g^{-1})$. Then $C_2\wr G$ is a group isomorphic to the semidirect product $C_2^Z\rtimes G$.\\

The representation theory of $C_2\wr G$ may be obtained by mean of the general representation theory of wreath products \cite{Hu,JK}, or, equivalently, by mean of the Frobenius-Mackey-Wigner theory of semidirect products with an abelian normal subgroup \cite{Serre, Simon}. We describe it briefly. The group $G$ acts on the dual group $\widehat{C_2^Z}=\{\chi_\theta:\theta\in C_2^Z\}$ by setting: $g\chi_\theta(\omega)=\chi_\theta(g^{-1}\omega)$, that is $g\chi_\theta=\chi_{g\theta}$.  The action of $G$ on $\widehat{C_2^Z}$ is equivalent to the action on $C_2^Z$ and both are the same thing as the action on the subsets of $Z$. In particular, the stabilizer
$G_\theta=\{g\in G:g\chi_\theta =\chi_\theta\}$ coincides with the stabilizer of $Z_\theta=\{z\in Z:\theta(z)=0\}$. The character $\chi_\theta$ has an {\em extension} to a character $\tilde{\chi}_\theta$ of $C_2\wr G_\theta$, defined by setting: $\tilde{\chi}_\theta(\omega,g)=\chi_\theta(\omega)$, for all $\omega\in C_2^Z$, $g\in G_\theta$. Similarly, if $\eta \in \widehat{G_\theta}$ (that is $\eta$ is an irreducible representation of $G_\theta$) then its {\em inflation} $\eta^\#$ to $C_2\wr G_\theta$ is defined by setting: $\eta^\#(\omega,g)=\eta(g)$, for all $\omega\in C_2^Z$, $g\in G_\theta$. Both $\tilde{\chi}_\theta$ and $\eta^\#$ are irreducible $C_2\wr G_{\theta}$-representations, and so is their tensor product $\tilde{\chi}_\theta\otimes \eta^\#$; clearly $\tilde{\chi}_\theta\otimes \eta^\#(\omega,g)=\chi_\theta(\omega)\eta(g)$. Now we can enunciate the main theorem in the representation theory of $C_2\wr G$.

\begin{theorem}\label{irrepwr}
Let $\Theta$ be a systems of representatives for the orbits of $G$
on $C_2^Z$ (any orbit has exactly one element in $\Theta$). Then

\[
\widehat{C_2\wr G}=\left\{\text{Ind}_{C_2\wr G_\theta}^{C_2\wr G}\tilde{\chi}_\theta\otimes \eta^\#: \quad\theta\in \Theta \quad\text{and}\quad \eta \in \widehat{G_\theta}\right\},
\]

that is the right hand side is a complete list of irreducible inequivalent representations of $C_2\wr G$.
\end{theorem}

Now suppose that $X$ is another homogeneous $G$-space. Fix $x_0\in X$ and set $H=\{g\in G:gx_0=x_0\}$, so that $X=G/H$. The group $C_2\wr G$ acts on $C_2^Z\times X$ by setting

\[
(\omega,g)(\theta,x)=((\omega,g)\theta,gx),\qquad\text{where}\qquad (\omega,g)\theta=\omega+g\theta
\]

for $(\omega,g)\in C_2\wr G$, $\theta\in C_2^Z$ and $x\in X$. We want to decompose the permutation representation of $C_2\wr G$ on $C_2^Z\times X$ into irreducible representations. Note that $L(C_2^Z\times X)\equiv L(C_2^Z)\otimes L(X)$. Moreover, the stabilizer of $(\mathbf{0}_Z,x_0)$ is just the subgroup $\tilde{H}=\{(\mathbf{0}_Z,h):h\in H\}\cong H$, that is $C_2^Z\times X\equiv (C_2\wr G)/\tilde{H}$. We begin with a general lemma on the action on a tensor product of the kind $\chi_\theta\otimes f$.

\begin{lemma}\label{genident}
If $(\omega,g)\in C_2\wr G$, $\theta\in C_2^Z$ and $f\in L(X)$ then

\[
(\omega,g)(\chi_\theta\otimes f)=\chi_{g\theta}(\omega)\cdot[\chi_{g\theta}\otimes gf].
\]
\end{lemma}
\begin{proof}
If $(\sigma,x)\in C_2^Z\times X$ then

\begin{equation*}
\begin{split}
[(\omega,g)(\chi_{\theta}\otimes f)](\sigma,x)=&(\chi_{\theta}\otimes f)[(\omega,g)^{-1}(\sigma,x)]\\
=&(\chi_{\theta}\otimes f)(g^{-1}\omega+g^{-1}\sigma,g^{-1}x)\\
=&\chi_{\theta}(g^{-1}\omega+g^{-1}\sigma)\cdot f(g^{-1}x)\\
=&\chi_{g\theta}(\omega)\cdot[\chi_{g\theta}\otimes gf](\sigma,x).
\end{split}
\end{equation*}
\end{proof}

For any $\theta \in \Theta$, choose a system $S_\theta$ of representatives for the left cosets of $G_\theta$ in $G$, that is $G=\coprod_{s\in S_\theta}sG_\theta$ (disjoint union). We always suppose that $1_G\in S_\theta$. For the moment, fix $\theta\in \Theta$ and suppose that $V$ is a $G_\theta$-invariant and irreducible subspace of $L(X)$.
We denote by $\eta$ the corresponding representation in $\widehat{G_\theta}$; but if $f\in V$ and $g\in G$ then the $g$-translate of $f$ is denoted by $gf$. Then the following corollary is an immediate consequence of Lemma \ref{genident}

\begin{corollary}\label{basicidentity}
If $(\omega,g)\in C_2\wr G$, $s\in S_\theta$, $gs=th$ with $h\in G_\theta$ and $t\in S_\theta$, and $f\in sV$ then

\[
(\omega,g)(\chi_{s\theta}\otimes f)=\chi_{t\theta}\otimes f'.
\]

where $f'=\chi_{gs\theta}(\omega)ths^{-1}f\in tV$.
\end{corollary}

\begin{lemma}\label{prodscal}
Suppose that $\theta'\in \Theta$, $s\in S_\theta$, $s'\in S_{\theta'}$ and that $V'$ is another $G_{\theta'}$-invariant subspace in $L(X)$. Then for $f\in sV$, $f'\in s'V'$ we have

\[
\langle\chi_{s\theta}\otimes sf,\chi_{s'\theta'}\otimes s'f' \rangle_{L(C_2^Z\times X)}=\delta_{\theta,\theta'}\delta_{s,s'}2^{\lvert Z\rvert}\langle f, f'\rangle_{L(X)}.
\]
\end{lemma}
\begin{proof}

We have $s\theta=s'\theta'$ if and only if $\theta=\theta'$ and $s=s'$. Therefore
\begin{equation*}
\begin{split}
\langle\chi_{s\theta}\otimes sf,\chi_{s'\theta'}\otimes s'f' \rangle_{L(C_2^Z\times X)}&=\langle \chi_{s\theta},\chi_{s'\theta'} \rangle_{L(C_2^Z)}\langle sf, s'f'\rangle_{L(X)}\\
&=\delta_{\theta,\theta'}\delta_{s,s'}2^{\lvert Z\rvert}\langle f, f'\rangle_{L(X)}.
\end{split}
\end{equation*}

\end{proof}

\begin{lemma}\label{indrep}
The space $\oplus_{s\in S_\theta}\{\chi_{s\theta}\otimes f:f\in sV\}$ is $C_2\wr G$-invariant and it is isomorphic to the irreducible representation $\text{Ind}_{C_2\wr G_\theta}^{C_2\wr G}\tilde{\chi}_\theta\otimes \eta^\#$.
\end{lemma}

\begin{proof}
From Corollary \ref{basicidentity} it follows that the subspace $\{\chi_{\theta}\otimes f:f\in V\}$ is $C_2\wr G_{\theta}$-invariant; moreover, the corresponding $C_2\wr G$-representation is equivalent to $\tilde{\chi}_\theta\otimes \eta^\#$. From the same corollary, it follows that the space $\oplus_{s\in S_\theta}\{\chi_{s\theta}\otimes f:f\in sV\}$ coincides with $\oplus_{s\in S_\theta}s\{\chi_{\theta}\otimes f:f\in V\}$ and that it is $C_2\wr G$-invariant. From Lemma \ref{prodscal} it follows that it is an orthogonal direct sum. Therefore we have verified all the requirements in the definition of induced representation \cite{Serre} (note also that $S_\theta$ is a system of representatives for the right cosets of $C_2\wr G_\theta$ in $C_2\wr G$).
\end{proof}

Now suppose that, for each $\theta\in\Theta$,

\begin{equation}\label{Xdec}
 L(X)=\bigoplus\limits_{i=0}^{n(\theta)}m_{\theta,i}V_{\theta,i}
\end{equation}

is the decomposition of $L(X)$ into irreducible $G_{\theta}$-representations. For different values of $i$ we have inequivalent representations and $m_{\theta,i}$ is the multiplicity of $V_{\theta,i}$ in $L(X)$. We also suppose that

\begin{equation}\label{Vdec}
m_{\theta,i}V_{\theta,i}=V^1_{\theta,i}\oplus V^2_{\theta,i}\oplus\cdots V^{m_{\theta,i}}_{\theta,i}
\end{equation}

is an explicit orthogonal decomposition of the isotypic block $m_{\theta,i}V_{\theta,i}$, (each $V^j_{\theta,i}$ is equivalent to $V_{\theta,i}$). For each $V^j_{\theta,i}$, set $W^j_{\theta,i}=\bigoplus\limits_{s\in S_\theta}\left\{\chi_{s\theta}\otimes f:f\in sV^j_{\theta,i}\right\}$. That is, $W^j_{\theta,i}$ is constructed as in Lemma \ref{indrep}, setting $V=V^j_{\theta,i}$. From Theorem \ref{irrepwr}, it follows that
all the representations $W^1_{\theta,i},W^2_{\theta,i},\dotsc W^{m_{\theta,i}}_{\theta,i}$ are irreducible and equivalent; by Lemma \ref{prodscal}, they are also mutually orthogonal subspaces of $L(C_2^Z\times X)$. We denote by $m_{\theta,i}W_{\theta,i}= W^1_{\theta,i}\oplus W^2_{\theta,i}\oplus\cdots W^{m_{\theta,i}}_{\theta,i}$ their direct sum.

\begin{theorem}\label{maintheorem}
The following
\begin{equation}\label{maindecomp}
L(C_2^Z\times X)=\bigoplus_{\theta\in\Theta}\bigoplus_{i=0}^{n(\theta)}m_{\theta,i}W_{\theta,i}
\end{equation}
is the decomposition of $L(C_2^Z\times X)$ into irreducible $C_2\wr
G$ representations and $m_{\theta,i}W_{\theta,i}=
W^1_{\theta,i}\oplus W^2_{\theta,i}\oplus\cdots
W^{m_{\theta,i}}_{\theta,i}$ is an orthogonal decomposition of the
isotypic block $m_{\theta,i}W_{\theta,i}$.
\end{theorem}

\begin{proof}
Another application of Lemma \ref{prodscal} yields the orthogonality of the decomposition \eqref{maindecomp}.
It remains only to show that the sum of all the spaces in right hand side of \eqref{maindecomp} is equal to $L(C_2^Z\times X)$. This is easy:

\begin{equation*}
\sum_{\theta\in\Theta}\sum_{i=0}^{n(\theta)} m_{\theta,i}\text{dim}W_{\theta,i}=
\sum_{\theta\in\Theta}\sum_{i=0}^{n(\theta)}\lvert S_{\theta}\rvert m_{\theta,i}\text{dim}V_{\theta,i}=\sum_{\theta\in\Theta}\left\lvert\frac{G}{G_\theta}\right\rvert\cdot\lvert X\rvert=2^{\lvert Z\rvert}\lvert X\rvert=\text{dim}L(C_2^Z\times X).
\end{equation*}
\end{proof}

The following corollary is a trivial consequence of Theorem \ref{maintheorem}, but it is worthwhile to enunciate it explicitly.

\begin{corollary}
The multiplicity of $\text{Ind}_{C_2\wr G_\theta}^{C_2\wr G}\tilde{\chi}_\theta\otimes \eta^\#$ in $L(C_2^Z\times X)$ is equal to the multiplicity of $\eta$ in the decomposition of $L(X)$ under the action of $G_\theta$.
\end{corollary}

Now we want to connect Theorem \ref{maintheorem} with the spectral analysis of an invariant operator.

\begin{proposition}
Let $\mathcal{M}:L(C_2^Z\times X)\rightarrow L(C_2^Z\times X)$ be a linear, selfadjoint, $C_2\wr G$-invariant operator.
\begin{enumerate}
\item
For any $\theta\in C_2^Z$, there exists a $G_\theta$-invariant, linear, selfadjoint operator
$M_\theta :L(X)\rightarrow L(X)$ such that:

\[
\mathcal{M}(\chi_\theta\otimes f)=\chi_\theta\otimes M_\theta f,
\]
for all $f\in L(X)$.
\item
Suppose that $V^j_{\theta,i}$ in \eqref{Vdec} is an eigenspace of $M_\theta$, with eigenvalue $\lambda^j_{\theta,i}$. Then the corresponding space $W^j_{\theta,i}$ in \eqref{maindecomp} is an eigenspace of $\mathcal{M}$, with the same eigenvalue $\lambda^j_{\theta,i}$.
\end{enumerate}
\end{proposition}

\begin{proof}
From Lemma \ref{genident} and the $C_2\wr G$-invariance of $\mathcal{M}$, we have:

\begin{equation}\label{actionprod}
(\omega,g)\mathcal{M}(\chi_\theta\otimes f)=\chi_{g\theta}(\omega)\cdot\mathcal{M}(\chi_{g\theta}\otimes gf).
\end{equation}

Setting $g=1_G$, \eqref{actionprod} becomes

\[(\omega,1_G)\mathcal{M}(\chi_\theta\otimes f)=\chi_\theta(\omega)\mathcal{M}(\chi_\theta\otimes f).
\]

This means that $\mathcal{M}(\chi_\theta\otimes f)$ belongs to the $\chi_\theta$-isotypic component in the decomposition of $L(C_2^Z\times X)$ under the action of $C_2^Z$, and therefore for any $f\in L(X)$ there exists
$f'\in L(X)$ such that:
$\mathcal{M}(\chi_{\theta}\otimes f)=\chi_\theta\otimes f'$. Setting $M_\theta f=f'$, we get a linear, selfadjoint operator $M_\theta :L(X)\rightarrow L(X)$ such that $\mathcal{M}(\chi_\theta\otimes f)=\chi_\theta\otimes M_\theta f$.

On the other hand, setting $\omega =\mathbf{0}_Z$ in \eqref{actionprod}, we get

\[
\chi_{g\theta}\otimes gM_\theta f=\chi_{g\theta}\otimes M_{g\theta}(gf),
\]

and therefore
$gM_\theta f=M_{g\theta}(gf)$. In particular, $M_\theta$ is $G_\theta$-invariant. Moreover, if $M_\theta f=\lambda^j_{\theta,i}f$ for all $f\in V^j_{\theta,i}$, then
also $M_{g\theta}(gf)=\lambda^j_{\theta,i}gf$.
From this fact it follows easily that $W^j_{\theta,i}$ is an eigenspace of $\mathcal{M}$, with the same eigenvalue $\lambda^j_{\theta,i}$.

\end{proof}

\begin{remark}
{\rm Clearly, the diagonalization of $M_{s\theta}$ is the same thing
as the diagonalization of $M_\theta$. If any $V_{\theta,i}^j$ in
\eqref{Xdec} as an eigenspace of $M_\theta$, then the action of the
group $C_2\wr G$ collects together all the eigespaces
$\{\chi_\theta\otimes f: f\in sV_{\theta,i}^j\}$ into a unique
eigenspace of $\mathcal{M}$, which is also an irreducible
representation.}
\end{remark}

We end this section with a general proposition of Harmonic Analysis
on spaces with multiplicity. We do not assume the previous notation.
Now $G$ is a finite group, $X$ a homogeneous $G$-space and
$L(X)=\oplus_{\rho\in J}m_\rho V_\rho$ is the decomposition of
$L(X)$ into irreducible $G$-representations; $m_\rho>0$ is the
multiplicity of the representation $\rho$. Denote by
$\text{Hom}_G(L(X),L(X))$ the {\it commutant} of $L(X)$, that is the
algebra of all operators $T:L(X)\rightarrow L(X)$ that commute with
the action of $G$:

\[
gTf=Tgf
\]

for any $g\in G,f\in L(X)$. Clearly any isotypic component $m_\rho V_\rho$ is $T$-invariant, for any $T\in\text{Hom}_G(L(X),L(X))$. The center of $\text{Hom}_G(L(X),L(X))$ is the subalgebra
$\{S\in\text{Hom}_G(L(X),L(X)):ST=TS\quad\text{for any}\quad T\in\text{Hom}_G(L(X),L(X))\}$.

\begin{proposition}\label{multspaces}
An operator $T\in \text{Hom}_G(L(X),L(X))$ belongs to the center of

$\text{Hom}_G(L(X),L(X))$ if and only if any isotypic component $m_\rho V_\rho$ is an eigenspace of $T$.
\end{proposition}

\begin{proof}
From Schur's lemma, we know that

\begin{equation}\label{isomcom}
\text{Hom}_G(L(X),L(X))\cong \oplus_{\rho\in J}M_{m_\rho,m_\rho}(\mathbb{C}),
\end{equation}

where $M_{m_\rho,m_\rho}(\mathbb{C})$ is the algebra of all $m_\rho\times m_\rho$ matrices over $\mathbb{C}$; see \cite{Sternberg}. Now we make the isomorphism \eqref{isomcom} more explicit. Suppose that $m_\rho V_\rho=V_\rho^1\otimes \dotsb\otimes V_\rho^{m_\rho}$ is an explicit orthogonal decomposition of $m_\rho V_\rho$ into $G$-irreducible representations. Using the Schur's lemma, we can introduce a basis $\{T_{i,j}^\rho:\rho\in J,\quad i,j=1,2,\dotsc,m_\rho\}$ for the commutant $\text{Hom}_G(L(X),L(X))$ with the following properties:

\[
\text{Ker}T^\rho_{i,j}=(V_\rho^j)^\bot,\qquad\qquad \text{Ran}T_{i,j}^\rho=V_\rho^i,\qquad\qquad\text{and}\quad T_{i,j}^\rho T^\rho_{j,k}=T^\rho_{i,k}.
\]

Then for any $T\in \text{Hom}_G(L(X),L(X))$, there exits a unique set of coefficients $\alpha_{i,j}^\rho$ such that $T=\sum_{\rho\in J}\sum_{i,j=1}^{m_\rho}\alpha_{i,j}^\rho T_{i,j}^\rho$, and the map

\[
T\mapsto \oplus_{\rho\in J}(a_{i,j}^\rho)_{i,j=1,\dotsc,m_\rho}
\]

is an explicit form of \eqref{isomcom}. Then the proposition is clear: $T$ is in the center of the commutant if and only if there exists $(\lambda_\rho)_{\rho\in J}$ such that

\[
\alpha_{i,j}^\rho=\delta_{i,j}\lambda_\rho.
\]
\end{proof}

If the operator $T$ is not in the center of the commutant, its diagonalization requires a suitable explicit decomposition of each isotypic component. This is the case of the lamplighter random walks considered in this paper. Another way to formulate and prove Proposition \ref{multspaces} is through the isomorphism between $\text{Hom}_G(L(X),L(X))$ and the convolution algebra of bi-$K$-invariant functions on $G$; see \cite{CST3}.

\section{The lamplighter on the circle revisited.}

Consider the wreath product $C_2\wr C_n$. The cyclic group $C_n$ will be written additively and it will be identified with $\mathbb{Z}/n\mathbb{Z}$. If $k\in C_n$ then we will think of $k$ as an integer representing $k+n\mathbb{Z}$. We will denote by $C_2^n$ the set of all maps $\theta:C_n\rightarrow C_2$.
If $k\in C_n$ and $\theta\in C_2^n$ then $k\theta(j)=\theta(j-k)$ and the group operation in $C_2\wr C_n=\{(\theta,k):\theta\in C_2^n,k\in C_n\}$ is:

\[
(\theta,k)(\omega,h)=(\theta+k\omega, k+h).
\]
Note that, in our notation, $s(k\omega)=(s+k)\omega$.
Any irreducible representation of $C_n$ is a one dimensional character of the form:
$e_k(h)=\exp\left(2\pi i\frac{hk}{n}\right)$, $h,k\in C_n$.

Think of $\theta\in C_2^n$ as a function $\theta:\mathbb{Z}\rightarrow C_2$ satisfying $\theta(k+n)=\theta(k)$ for any $k\in \mathbb{Z}$. Then the {\em period} of $\theta$ is the smallest positive integer $t=t(\theta)$ such that $\theta(k+t)=\theta(k)$ for any $k\in\mathbb{Z}$; clearly $t$ divides $n$ and if $n=mt$ then the stabilizer of $\theta$ is the subgroup $C_m=\langle t \rangle$ (recall also that for any divisor $m$ of $n$, the subgroup of $C_n$ isomorphic to $C_m$ is unique \cite{Lang}). The characters of the subgroup $\langle t\rangle$ are given by the restrictions: $e_0|_{\langle t\rangle}, e_1|_{\langle t\rangle},\dotsc,e_{m-1}|_{\langle t\rangle}$, where $e_0, e_1,\dotsc,e_{m-1}$ are as above. Indeed, for $0\leq r,l\leq m-1$ we have: $e_r(lt)=\exp\left(2\pi i\frac{rlt}{n}\right)=\exp\left(2\pi i\frac{rl}{m}\right)$. We set $e_r|_{\langle t\rangle}(k) =e_r(k)$ when $k\in \langle t\rangle$, $e_r|_{\langle t\rangle}(k)=0$ otherwise.
In what follows, we also set $m(\theta)=\frac{n}{t(\theta)}$, but  we will write simply $t$ and $m$ when it is clear the $\theta$ we are talking about.

Now take $\theta\in C_2^n$ and $0\leq r\leq m-1$. If we compute the inflation of $e_r|_{\langle t\rangle}$ and the extension of $\chi_\theta$, we get the character $\tilde{\chi}_\theta\otimes (e_r|_{\langle t\rangle})^\#$ of $C_2^n\wr \langle t\rangle$ given by: $\tilde{\chi}_\theta\otimes (e_r|_{\langle t\rangle})^\#(\omega,lt)=\chi_\theta(\omega)e_r(lt)$, for $\omega\in C_2^n$ and $l=0,1,\dotsc,m-1$. Let $\Theta$ be a  set of representatives for the orbits of $C_n$ on $C_2^n$ (such orbits may be enumerated by mean of the so called Polya-Redfield theory; see \cite{Lint-Wi} for an elementary account and \cite{Kerber} for a more comprehensive treatment). Then we can apply Theorem \ref{irrepwr}.

\begin{theorem}
The set $\{\text{Ind}_{C_2\wr \langle t(\theta)\rangle}^{C_2\wr C_n}\left[\tilde{\chi}_\theta\otimes (e_r|_{\langle t(\theta)\rangle})^\#   \right]:\theta\in\Theta, r=0,1,\dotsc,m(\theta)-1\}$ is a complete list of irreducible inequivalent representations of $C_2^n\wr C_n$.
\end{theorem}

Suppose again that $\theta\in \Theta$. For $s=0,1,\dotsc,t-1$, set $\Omega_s=\{s,s+t,\dotsc,s+(m-1)t\}$, and for $r=0,1,\dotsc, m-1$,

\[
f_{r,s}(k)=se_{m-r}|_{\langle t\rangle}\equiv\left\{\begin{array}{lll}e_{m-r}(k-s)&\text{if}&k\in\Omega_s\\
0&\text{if}&k\notin \Omega_s.\end{array}\right.
\]

Clearly $f_{r,s}\in L(\Omega_s)$ and

\begin{equation}\label{hfrs}
hf_{r,s}=\left\{\begin{array}{lll}f_{r,s+h}&\text{if}&h\notin\langle t\rangle\\
e_r(h)f_{r,s}&\text{if}&h\in\langle t\rangle.\end{array}\right.
\end{equation}

But $C_n=\coprod_{s=0}^{t-1}\Omega_s$ is the decomposition of $C_n$ into $\langle t\rangle$-orbits, and therefore from \eqref{hfrs} it follows that

\[
L(C_n)= \bigoplus_{r=0}^{m-1}\langle f_{r,0},f_{r,1},\dotsc,f_{r,t-1}\rangle
\]

is the decomposition of $L(C_n)$ into $C_r$-isotypic components, where the $r$-th summand is precisely the $e_r|_{\langle t\rangle}$-isotypic component. Now consider the operator $\mathcal{M}$ of edge lamplighter random walk, as in Section \ref{disccircl}.
Clearly, $\langle f_{r,0} \rangle\oplus\langle f_{r,1}\rangle\oplus\dotsb\oplus\langle f_{r,t-1}\rangle$ is an orthogonal decomposition into irreducible representations, but $\mathcal{M}$ is {\em not} diagonal in this decomposition. Now we show that we need another application of \eqref{circlespec} and \eqref{specpath}.\\

In the notation of Section \ref{disccircl}, we can think of $\theta$ as a function defined on the vertices, by setting $\theta(k)=\theta(\{k,k+1\})$. Moreover, in the notation of the present section, we can always suppose that, for any $\theta\in\Theta$ with $\theta\neq\mathbf{0}_{C_n}$, we have $\theta(-1)=1$. Then the spectrum of $M_\theta$ is clearly $m$ times the spectrum of its restriction to $L(\{0,1,\dotsc,t-1\})$. Similarly, if $\alpha_s,s=0,1,\dotsc,t-1$ are complex numbers and $\alpha_{s+lt}=\alpha_s$ for any $l\in\mathbb{Z}$, an application of \eqref{hfrs} yields

\[
M_\theta(\alpha_0f_{r,0}+\alpha_1 f_{r,1}+\dotsb+\alpha_{t-1}f_{r,t-1})=\sum_{s=0}^{t-1}
\left(\frac{1-\theta(s-1)}{2}\alpha_{s-1}+\frac{1-\theta(s-1)}{2}\alpha_s\right)f_{r,s}.
\]

In other words, the eigenvalue problem of
$M_\theta|_{\langle f_{r,0},\dotsc,f_{r,t-1}\rangle}$ (with respect to the basis $\{ f_{r,0},\dotsc,f_{r,t-1}\}$)
coincides with the eigenvalue problem of $M_\theta|_{ L(\{0,1,\dotsc,t-1\})}$
(with respect to the basis $\{\delta_0,\delta_1\dotsc,\delta_{t-1}\}$).
Using \eqref{specpath} we can obtain an orthogonal decomposition $\langle f_{r,0},f_{r,1},\dotsc,f_{r,t-1} \rangle=\oplus_{j=0}^{t-1}\langle \phi_{r,j}\rangle$ such that any $\phi_{r,j}$ is an eigenvector of $M_\theta$. In the following proposition, we give the obvious conclusions of the preceding discussion.

\begin{proposition}
Suppose that $\phi_{r,0},\phi_{r,1},\dotsc,\phi_{r,t-1}$ are as above. Then

\[
\bigoplus_{j=0}^{t-1}\text{Ind}_{C_2\wr\langle t\rangle}^{C_2\wr C_n}\langle \tilde{\chi}_\theta \otimes (\phi_{r,j})^\# \rangle
\]

is a decomposition of the $\text{Ind}_{C_2\wr \langle t(\theta)\rangle}^{C_2\wr C_n}\left[\tilde{\chi}_\theta\otimes (e_r|_{\langle t(\theta)\rangle})^\#   \right]$-isotypic component of $L(C_2\wr C_n)$ into eigenspaces of the lamplighter operator $\mathcal{M}$.
\end{proposition}

Note that, for a fixed $\theta$, the eigenvalues do not depend on $r\in\{0,1,\dotsc,m\}$. Moreover, from Proposition \ref{multspaces} we deduce that $\mathcal{M}$ is not in the center of the group algebra of $C_2\wr C_n$.

\section{On a decomposition of Schoolfield.}
In this subsection, we want to apply Theorem \ref{maintheorem} to get the decomposition of the homogeneous space of the signed Bernoulli-Laplace contained in \cite{Schoolfield}. First of all, we need a description of the irreducible representations of the hyperoctahedral group $C_2\wr S_n$.
See also \cite{GK,JK}. Now $G=S_n$ and $Z=\{1,2,\dotsc,n\}$. For any $0\leq k\leq n$, choose $\theta^{(k)}\in C_2^Z$ such that $\lvert \{j\in Z:\theta^{(k)}(j)=0\}\rvert=k$. Then $\{\theta^{(0)},\theta^{(1)},\dotsc,\theta^{(n)}\}$ is a set of representatives for the orbits of $S_n$ on $C_2^Z$. Moreover, the stabilizer of $\theta^{(k)}$ is isomorphic to $S_k\times S_{n-k}$. We recall that the irreducible representations of the symmetric group $S_t$ are canonically parametrized by the partitions of $t$; \cite{JK,Sagan}. For $\lambda\vdash t$ (this means that $\lambda$ is a partition of $t$), we will denote by $\rho_\lambda$ the irreducible representation of $S_t$ canonically associated to $\lambda$ and by $S^\lambda$ the corresponding representation space.
As usual \cite{Sagan}, we set $M^{n-m,m}=L(S_n/(S_m\times S_{n-m}))$, that is $M^{n-m,m}$ is the permutation representation of $S_n$ on the space of all $m$-subsets of $\{1,2,\cdots,n\}$. We recall that

\begin{equation}\label{Mdec}
M^{n-m,m}=\bigoplus_{k=0}^{\min\{m,n-m\}}S^{n-k,k}.
\end{equation}

See \cite{JK}; see also \cite{CST1} for an elementary exposition.\\

We will use the following notations: if $A$ is a set with $\lvert A\rvert=k$ and $0\leq l\leq k$ then $M^{k-l,l}(A)$ will denote the space $M^{k-l,l}$ constructed by using the $l$-subsets of $A$ and

\begin{equation}\label{Mdec2}
M^{k-l,l}(A)=\bigoplus_{j=0}^{\min\{l,k-l\}}S_l^{k-j,j}(A).
\end{equation}

the corresponding decomposition into irreducible $S_k$-representations, as in \eqref{Mdec}. That is, $S_l^{k-j,j}(A)$ is the subspace of $M^{k-l,l}(A)$ isomorphic to $S^{k-j,j}$. In \cite{CST1}, the decomposition \eqref{Mdec2} is realized concretely by mean of finite Radon transforms; see also \cite{CST3,Diaconis,Sc-To2}.\\

The irreducible representations of the group $S_k\times S_{n-k}$ are all of the form $\rho_\lambda\otimes \rho_\mu$, for $\lambda\vdash k$ and $\mu\vdash n-k$.
If we set $\rho_{[\lambda;\mu]}=\text{Ind}_{C_2\wr(S_k\times S_{n-k})}^{C_2\wr S_n}[\tilde{\chi}_{\theta^{(k)}}\otimes (\rho_\lambda \otimes \rho_\mu)^\#]$, applying Theorem \ref{irrepwr} we can say that

\[
\{\rho_{[\lambda;\mu]}:\lambda\vdash k, \mu\vdash n-k\quad\text{and}\quad 0\leq k\leq n\}
\]

\noindent
is a complete list of inequivalent, irreducible $C_2\wr S_n$-representations.\\

Now fix $1\leq r \leq n-1$ and suppose that $X$ is the family of all $r$-subsets of $Z=\{1,2,\dotsc,n\}$. The homogeneous space of the signed Bernoulli-Laplace diffusion model studied in \cite{Schoolfield} coincides with $C_2^Z\times X$. Now we give a decomposition of the space $L(C_2^Z\times X)$ into irreducible $C_2\wr S_n$-representations.

\begin{theorem}
A decomposition of the permutation representation of $C_2\wr S_n$ on $C_2^Z\times X$ is given by:

\[
L(C_2^Z\times X)=
\bigoplus_{k=0}^n\bigoplus_{i=\max\{0,r+k-n\}}^{\min\{k,r\}}\bigoplus_{l=0}^{\min\{i,k-i\}}\bigoplus_{m=0}^{\min\{n-k-r+i,r-i\}}W^i_{k;l,m}
\]

where

\[
W^i_{k;l,m}=\langle \chi_\theta\otimes (f_1\otimes f_2)\in L(C_2^Z)\otimes L(X):\lvert Z_\theta\rvert=k,f_1\in S_i^{k-l,l}(Z_\theta)\quad\text{and}\quad f_2\in S_{r-i}^{n-k-m,m}(Z\setminus Z_\theta)\rangle.
\]

Moreover, the representation of $C_2\wr G$ on $W^i_{k;l,m}$ is isomorphic to
$\rho_{[(k-l,l);(n-k-m,m)]}$.

\end{theorem}
\begin{proof}
In order to apply Theorem \ref{maintheorem}, we need to decompose
the space $L(X)$ into irreducible $S_k\times
S_{n-k}$-representations, for any $0\leq k\leq n$. Suppose that
$B_k$ is the $k$-subset of $Z$ fixed by $S_k\times S_{n-k}$. Then
the orbits of $S_k\times S_{n-k}$ on $X$ are $\Xi_i=\{A\in X:\lvert
A\cap B_k\rvert=i\}$, that is the orbit of an element $A\in X_r$
(which is an $r$-subset of $Z$) is determined by the cardinality of
its intersection with $B_k$. Clearly $\max\{r+k-n,0\}\leq
i\leq\min\{k,r\}$ (intersect $A$ with $B_k$ and with the complement
of $B_k$). Applying \eqref{Mdec2}, we get:

\begin{equation*}
\begin{split}
L(\Xi_i)&\cong M^{k-i,i}(B_k)\otimes M^{n-k-r+i,r-i}(Z\setminus B_k)\\
&=\bigoplus_{l=0}^{\min\{i,k-i\}}\bigoplus_{m=0}^{\min\{n-k-r+i,r-i\}}
\left(S_i^{k-l,l}(B_k)\otimes S_{r-i}^{n-k-m,m}(Z\setminus B_k)\right).
\end{split}
\end{equation*}

Therefore, the permutation representation of $S_k\times S_{n-k}$ on $L(X)$ decomposes as follows:

\[
L(X)=\bigoplus_{i=\max\{0,r+k-n\}}^{\min\{k,r\}}
\bigoplus_{l=0}^{\min\{i,k-i\}}\bigoplus_{m=0}^{\min\{n-k-r+i,r-i\}}\left(S_i^{k-l,l}(B_k)\otimes S_{r-i}^{n-k-m,m}(Z\setminus B_k)\right)
\]

and an application of Theorem \ref{maintheorem} ends the proof.

\end{proof}

Just set $j=k-i$ to get exactly the formula of lemma 3.2.1 in \cite{Schoolfield}.
Summing up the equivalent representations, we can say that the decomposition of the permutation representation of $C_2\wr G$ on $C_2^Z\times X$ is given by:

\[
\bigoplus_{k=0}^n\bigoplus_{l=0}^{\min\{k,r,n-r\}}\bigoplus_{m=0}^{\min\{n-k,n-r-l,r-l\}}
m_{k;l,m}\rho_{[(k-l,l);(n-k-m,m)]}
\]

where $m_{k;l,m}=\min\{k-l,r-m\}-\max\{r+k-n+m,l\}$.

\section{The lamplighter on the complete graph revisited.}
Setting $r=1$ in the results of the previous subsection, we get an explicit decomposition for the vertex lamplighter on the complete graph. Now $X=Z$; moreover, $S_1^{(k)}(A)$ are the constant functions, while $S_1^{k-1,1}(A)$ is made up of the functions on $A$ satisfying $\sum_{a\in A}f(a)=0$. Therefore

\begin{equation*}
\begin{split}
&W^1_{k;0,0}=\langle\chi_\theta\otimes f: \quad\lvert X_\theta \rvert=k,\ f|_{X_\theta}\in S_1^{(k)}(X_\theta) \quad\text{and}\quad f|_{X\setminus X_\theta}\equiv 0\rangle\\
&W^0_{k;0,0}=\langle\chi_\theta\otimes f:\quad\lvert X_\theta \rvert=k, \quad
f|_{X_\theta}\equiv 0\quad\text{and}\quad f|_{X\setminus X_\theta}\in S_1^{(n-k)}(X\setminus X_\theta)\rangle,\\
&W^1_{k;1,0}=\langle\chi_\theta \otimes f:\quad\lvert X_\theta \rvert=k, f|_{X_\theta}\in S_1^{k-1,1}(X_\theta)\quad\text{and}\quad f|_{X\setminus X_\theta}\equiv 0\rangle,\\
&W^0_{k;0,1}=\langle \chi_\theta\otimes f: \quad\lvert X_\theta \rvert=k,
f|_{X_\theta}\equiv 0
\quad\text{and}\quad f|_{X\setminus X_\theta}\in S_1^{n-k-1,1}(X\setminus X_\theta)\rangle.
\end{split}
\end{equation*}

Then we can write the decomposition of $L(C_2^X\times X)$ into irreducible $C_2\wr S_n$-representations:

\[
L(C_2^X\times X)=\left(W^0_{0;0,0}\oplus W^0_{0;0,1}\right)\oplus\left[\bigoplus_{k=1}^{n-1}\left(W^0_{k;0,0}\oplus W^1_{k;0,0}\oplus W^0_{k;0,1}\oplus W^1_{k;1,0}\right)\right]\oplus \left(W^1_{n;0,0}\oplus W^1_{n;1,0}\right).
\]

In particular, the representations $W_{k;1,0}$ and $W_{k;0,1}$ have multiplicity 1, while the representations $W_{k;0,0}$ have multiplicity 2. Moreover, in the notation of Section \ref{completegraph} we have:

\[
W^1_{k;0,0}=\bigoplus_{\substack{\theta\in C_2^X\\\lvert X_\theta\rvert=k}}\text{Ran}(P_\theta),\qquad\qquad W^1_{k;1,0}=\bigoplus_{\substack{\theta\in C_2^X:\\\lvert X_\theta\rvert=k}}\text{Ran}(R_\theta-P_\theta)
\]

while

\[
\left(
\bigoplus_{k=0}^{n-1}W^0_{k;0,0}
\right)\bigoplus\left(\bigoplus_{k=1}^{n-1}W^0_{k;0,1}\right)
\]

is the decomposition of the null eigenspace into irreducible representations. Again, the operator is not in the center of the commutant algebra: $W^0_{k;0,0}$ and $W^1_{k;0,0}$ are equivalent but they correspond to different eigenvalues, namely $0$ and $\frac{k-1}{n-1}$.

\begin{remark}\label{lastrem}{\rm
Consider the following mixing procedure for the lamplighter on the
complete graph: at each time a random pair of distinct vertices
$x,y$ is chosen. Both the lamps in $x$ and $y$ are randomized.
Moreover, if the lamplighter is in $x$ (resp. $y$), it moves to $y$
(resp. $x$); if the lamplighter is in $X\setminus\{x,y\}$, then it
remains in his position. This is a slight variation of the mixing
procedure in \cite{Schoolfield}, for $r=1$. Now the corresponding
Markov operator $\mathcal{M}'$ is in the center of the commutant
algebra: it is easy to show that $\mathcal{M}'(\chi_\theta\otimes
f)=\chi_\theta\otimes M'_\theta f$, where

\[
M'_\theta f(x)=\frac{2}{n}M_\theta f(x)+\frac{\lvert X_\theta\setminus\{x\}\rvert(\lvert X_\theta\setminus\{x\}\rvert-1)}{n(n-1)}f(x),
\]

and that the whole $W^0_{k;0,0}\oplus W^1_{k;0,0}$ is an eigenspace, with corresponding eigenvalue equal to $\frac{k(k-1)}{n(n-1)}$. Define $\mathcal{C}$ as the set of all pairs $(\theta,\tau)\in C_2\wr S_n$ such that $\tau$ is a transposition and $\theta(x)=0$ if $\tau(x)=x$. Then $\mathcal{C}$ is a conjugacy class of $C_2\wr S_n$ \cite{JK} and  $\mathcal{M}'F(\omega,x)=\frac{1}{2n(n-1)}\sum_{(\theta,\tau)\in\mathcal{C}}F(\omega+\theta,\tau (x))$.
This is the reason for which $\mathcal{M}'$ is in the center of the commutant.
}
\end{remark}

\qquad\\
\qquad\\
\noindent
FABIO SCARABOTTI, Dipartimento MeMoMat, Universit\`a li Studi di Roma ``La Sapienza'', via A. Scarpa 8, 00161 Roma (Italy)\\
{\it e-mail:} {\tt scarabot@dmmm.uniroma1.it}\\
FILIPPO TOLLI, Dipartimento di Matematica, Universit\`a Roma TRE, L. San Leonardo Murialdo 1, 00146 Roma, Italy
{\it e-mail:} {\tt tolli@mat.uniroma3.it}\\

\end{document}